\documentclass[12pt,twoside,english,shortlabels]{amsart}

\usepackage{babel}
\usepackage{amstext}
\usepackage{amssymb}

\usepackage{amsmath} 
\usepackage{latexsym}
\usepackage{graphicx} 

\usepackage{mathptmx}

\usepackage{graphicx} 

\makeatletter

\usepackage{amsthm}
\usepackage{enumitem}		% customizable list environments
\usepackage{mathrsfs} %for the \mathscr
\providecommand{\noopsort}[1]{}

\numberwithin{equation}{subsection}

\theoremstyle{definition} 

 \newtheorem{definition}{Definition}[section]
 \newtheorem{remark}[definition]{Remark}
 \newtheorem{example}[definition]{Example}
 
 \newtheorem{Conjecture}[definition]{Conjecture}

\newtheorem*{notation}{Notations}

\theoremstyle{plain}      

 \newtheorem{proposition}[definition]{Proposition}
 \newtheorem{theorem}[definition]{Theorem}

\newtheorem{conjecture}{Conjecture}

\makeatletter
\newcommand*{\house}[1]{
  \mathord{
    \mathpalette\@house{#1}
  }
}
\newcommand*{\@house}[2]{
  \dimen@=\fontdimen8 %
      \ifx#1\scriptscriptstyle\scriptscriptfont
      \else\ifx#1\scriptstyle\scriptfont
      \else\textfont\fi\fi
      3 %
  \sbox0{%
    $#1%
      \vrule width\dimen@\relax
      \overline{%
        \kern2\dimen@
        \begingroup % to keep changes of \dimen@ in #2 local
          #2%
        \endgroup
        \kern2\dimen@
      }
      \vrule width\dimen@\relax
      \mathsurround=1.5\dimen@ % outside margin
    $
  }
  \ht0=\dimexpr\ht0-\dimen@\relax
  \dp0=\dimexpr\dp0+2\dimen@\relax
  \vbox{
    \kern\dimen@ 
    \copy0 
  }
}

\def\dyg{{\rm dyg}}

\newcommand{\nb}{\mathbb{N}}

\let\cal\mathcal

\def\11{{\mathbf 1}}

\theoremstyle{remark}
\newtheorem{exampl}[subsubsection]{Example}

\def\bee{\begin{exampl}}
\def\eee{\end{exampl}}
\def\bn{\begin{notation}}
\def\en{\end{notation}}
\def\br{\begin{remark}}
\def\er{\end{remark}}
\def\bp{\begin{prop}}
\def\ep{\end{prop}}
\def\bpr{\begin{proof}}
\def\epr{\end{proof}}
\def\bt{\begin{thm}}
\def\et{\end{thm}}
\def\be{\begin{equation}}
\def\ee{\end{equation}}
\def\bl{\begin{lem}}
\def\el{\end{lem}}
\def\bc{\begin{cor}}
\def\ec{\end{cor}}
\def\bd{\begin{defn}}
\def\ed{\end{defn}}

\numberwithin{equation}{subsection}

\author{Jean-Louis Verger-Gaugry}
\thanks{}
\address{Univ. Grenoble Alpes, Univ. Savoie Mont-Blanc, CNRS UMR 5127, LAMA, 
F-73000 Chamb\'ery, France}
\email{verger-gaugry.jean-louis@unsechat-math.fr}

\begin{document}

\dedicatory{\today}

%-----------------------------------------------------
%%%%%%%%%%%%%%%%%%%%%JLVG%%%%%%%%%%%%%%%%%%%%%%

\def\dyg{{\rm dyg}}

%%%%%%%%%%%%%%%%%%%%%%%%%%%%%%%%%%%%%%%%%%%%%%%

\title{Divergent Trajectories in the $3x+1$ Problem
and a $S$-units Equation}

\begin{abstract}
In this note we assume the existence of divergent
trajectories (if any) under the Collatz map
and investigate the consequences 
of this assumption, even though
the Conjecture of Collatz seems to be true.
In a trajectory we show that the Diophantine 
exponential equation 
$x+y=w 2^n$ in the integers
$x,y,w,n$ plays an important role.
In this context the recent strong finiteness
theorem of Bennett and Billerey 
on $S$-units has consequences on
 the set of prime numbers
generated by a collection of
local minima, called wave-minima
(defined ad hoc for simplicity's sake)
in a divergent trajectory.
It allows to define canonically
a density of a trajectory and some 
conjectures and open problems
related to these prime numbers.

\end{abstract}

\maketitle

\vspace{0.2cm}


\vspace{0.4cm}


\bigskip

\section{Introduction}
\label{S1}

The Collatz function
$T : \nb \to \nb$ is defined by
\begin{equation}\label{eq T}
T(n) = 
\begin{cases} 
n/2 & \textrm{ if } n \equiv 0 \bmod 2, \\
(3n+1)/2 & \textrm{ if } n \equiv 1 \bmod 2.
\end{cases}
\end{equation}
This function became famous due to the 
Syracuse Conjecture, also named
$3x+1$ Conjecture or Collatz Conjecture.
This Conjecture has been formulated initially by
Collatz as a numerical observation dating back to
the 1930's \cite{chamberland}.
Let us recall it. Denote
$T^{(0)}:=$Id, 
$T^{(1)}:=T$,
$T^{(j)}:= T\circ T^{(j-1)}$
for $j \geq 2$.

\begin{conjecture}
Let $n \in \nb$, $n \geq 1$. Then there exists
$k \geq 1$ such that 
$T^{(k)}(n)=1$.
\end{conjecture}
This Conjecture means that any trajectory
$(T^{(j)}(n))_{j \geq 0}$, $n \geq 1$, is
bounded 
(by an upper bound depending upon $n$)
and ends in the trivial cycle
$\{1,2\}$. We refer the reader to
\cite{allouche}
\cite{chamberland} \cite{lagarias}
\cite{eliahou}
for more details about the Conjecture.

The Conjecture has been recently
verified on a large
number of integers $n$ by Barina \cite{barina}.
Few is known on the boundedness problem of the
trajectories, though Tao 
recently brought some light on it in \cite{tao}.
Whether divergent trajectories exist or not is 
not known. Proving that divergent
trajectories do not exist would be a 
formidable result.

\

In this note we proceed a contrario.
We assume that a divergent trajectory
exists and investigate the consequences
of this assumption.
In this divergent trajectory we show that the 
Diophantine exponential equation
\begin{equation}
\label{diophantexp}
x+y=w 2 ^n
\end{equation}
with $x, y, w, n$ integers plays an important role
in the successive iterations of the Collatz function
$T$, in strings of successive odd integers
of height greater than $2$ 
(the height is defined below).
We show that this role is due to a recent 
finiteness
theorem obtained by Bennett and Billerey.
This theorem, Theorem 3.3 in \cite{bennettbillerey},
is recalled below as
Theorem \ref{thmBennettBillerey},
in which the equation 
\eqref{diophantexp} is an $S$-units
equation, i.e. in which the solutions 
$x,y,w,n$ have prime decompositions 
for which the prime factors are all in a
 finite set of primes $S$.
Finiteness means that these authors, 
using strong theorems, prove that
a subset of solutions
of this equation is necessarily finite
(cf Section \ref{S2} for the details).

In this note we report the rather suprising
discovery that
this theorem of Bennett and Billerey has consequences
on the prime numbers which appear as prime factors
in the prime decompositions of the iterates of $T$
which constitute the
strings of successive odd integers in the 
divergent trajectory. 
We notice that the multiplicities of these
prime factors occurring in these prime decompositions
do not play any role.
The links
between Bennett-Billerey's Theorem and
the concerned prime numbers are explicited
in Section \ref{S2}.

As a consequence, in 
Section \ref{S3}, a canonical definition of
the density of a divergent trajectory, using 
these prime factors, is
stated. It leads to state a Conjecture
on divergent trajectories (if any) which 
have density zero. Some open problems 
are formulated in
Section \ref{S3}.

\

For doing this, two types of local minima of a 
given trajectory are defined. These definitions
are chosen by the author ad hoc and for 
simplicity's sake.
Let us fix the notations.

\

Let $n \geq 3$ be on odd integer. 
Let $\mathcal{T}:= (T^{(j)}(n))_{j \geq 0}$
denote the trajectory of $n$. By
definition a local minimum
of $\mathcal{T}$ is either 
$n$ itself or an integer
$T^{(j)}(n)$ for some $j \geq 1$
such that
$$T^{(j-1)}(n) > T^{(j)}(n) < T^{(j+1)}(n).$$
A local minimum is odd.
If $T^{(j)}(n)$, $j \geq 0$, is a local minimum,
the iterate $T^{(j+1)}(n)=T \circ T^{(j)}(n)$ 
is either even or odd. By definition
a {\em ripple} of $\mathcal{T}$ is 
a pair $(T^{(j)}(n), 
T^{(j+1)}(n))$ 
in which
$T^{(j)}(n)$ is a local minimum and
$T^{(j+1)}(n)$ is even.
By definition,
a {\em wave of height $k \geq 2$}
is a $(k+1)$-tuple of integers
$(T^{(j)}(n), T^{(j+1)}(n), \ldots, T^{(j+k)}(n))$
in which $T^{(j)}(n)$ is a local minimum
of $\mathcal{T}$ and for which
$T^{(j+s)}(n)$, for $s=1$ to $k-1$, are all
odd with $T^{(j+k)}(n)$ even.
By {\em wave} we mean any wave of height 
$k$, for some integer $k \geq 2$.
Any local minimum $\mu$ of
$\mathcal{T}$ is either the 
first element of a 
ripple or of a wave: 
in the first case
we say that $\mu$ is a 
{\em ripple-minimum}, in the
second case, that $\mu$ is
a {\em wave-minimum}.

\

Let $\mathbb{P}$ denote the set of prime numbers.
Let $\mu$ be a wave-minimum and denote by
$(\mu = T^{(j)}(n), 
T^{(j+1)}(n), \ldots, T^{(j+k)}(n))$ the 
wave whose first element is $\mu$.
By a prime number generated by the wave-minimum 
$\mu$
we mean a prime number occurring either
in the prime decomposition of $\mu+1$,
or in the prime decomposition of 
$T^{(j+s)}(n)$, for all $s$ such that $1 \leq s \leq k$ . If
$\mu + 1 = p_{1}^{\alpha_1}
p_{2}^{\alpha_2}\ldots
p_{q}^{\alpha_q}$, we denote:
$${\mathcal Prm}(\mu)=
\{p_1, p_2 , \ldots, p_q\}~
\bigcup
~~\cup_{s=1}^{k}
\{q \in \mathbb{P} : q 
~~\mbox{is a prime factor of}
~~T^{(j+s)}(n)\} ~\bigcup ~\{3\}$$
the set of the prime numbers
generated by $\mu$.

\

Let us examplify the terminology introduced
on some non-divergent trajectories.

\begin{example}
a) The trajectory of $n=17$ under $T$ 
is the following:
\begin{equation}
\label{traj17}
17 \to 26 \to 13 \to 20 \to 10 \to 5 \to 8 \to 4 \to 2 \to 1 \to 2 \to 1 \to 2 \to \ldots.
\end{equation}
It has 3 ripples, successively: $\{T^{(0)}(17)=17, 
T^{(1)}(17)=26\}$,
$\{T^{(2)}(17)=13, T^{(3)}(17)=20\}$, 
$\{T^{(5)}(17)=5, T^{(6)}(17)=8\}$,
and 
3 ripple-minima greater than 1: 
$\{17\}, \{13\}, \{5\}$. The trajectory
of $17$ has no wave and no wave-minimum
and infinitely many 
ripples $\{1,2\}$.

\

b) The trajectory of $n=19$ under $T$ contains
waves and wave-minima.
It is
$$19 \to 29 \to 44 \to 22 \to 11 \to 17 \to 26 \to 
\ldots$$
ending in the infinite sequence
\eqref{traj17} of the iterates of 17.
The two waves of the trajectory of $19$, 
each of them of height 2, are:
$\{19, 29, 44\}$ and
$\{11, 17, 26\}$. The two wave-minima are
$\{\mu_{1,w}=19\}, \{\mu_{2,w}=11\}$. Since
$19+1= 20 = 2^2 \times 5$
and $11+1=12=2^2 \times 3$, 
we have:
$${\cal Prm}(\mu_{1,w})=\{2, 5\}
\cup
\{29\}
\cup
\{2, 11\} \cup \{3\}$$
\hspace{0.5cm}
and
\hspace{0.5cm}
$${\cal Prm}(\mu_{2,w})=\{2, 3\}
\cup
\{17\}
\cup
\{2, 13\}.$$
\end{example}

\

Let 
$\mathcal{T}=
(T^{(j)}(n))_{j \geq 0}$ be a divergent
trajectory for some odd integer
$n \geq 3$. The wave-minima of $\mathcal{T}$
uniquely define the infinite increasing 
sequence of integers
$\mathcal{J}= \{i_j : j= 1, 2, 3, \ldots\})$,
$0 \leq i_1 < i_2 <i_3 < \ldots,$
for which the
$T^{(i_j)}(n)$s are exactly the successive
wave-minima
of $\mathcal{T}$.
We enumerate the successive wave-minima of
$\mathcal{T}$ as
follows:
$$\mu_{1,w} := T^{(i_1)}(n),~~
\mu_{2,w} := T^{(i_2)}(n),~~
\mu_{3,w} := T^{(i_3)}(n),
\ldots.$$

\

If $\kappa$ denotes a positive integer
and $p_{1}^{\beta_1} p_{2}^{\beta_2} 
\ldots p_{q}^{\beta_q}$ its prime decomposition,
we  call
$\beta_i$  the $p_i$-adic valuation of
$\kappa$, and denote it by 
$v_{p_i}(\kappa)$.

\

\begin{theorem}
\label{mainthm}
Let 
$\mathcal{T}=
(T^{(j)}(n))_{j \geq 0}$ be a divergent
trajectory for some odd integer
$n \geq 3$.
Then
\begin{enumerate}
\item[(i)] the ripples of $\mathcal{T}$, possibly 
in infinite number, do not
contribute to the divergence,
\item[(ii)] any wave 
$(\mu_{j,w}=T^{(i_j)}(n),
T^{(i_j + 1)}(n),\ldots,
T^{(i_j +k)}(n))$, for some $j \geq 1$
and $k \geq 2$,
has a finite height 
$k = v_{2}(\mu_{j,w}+1)$ which satisfies
\begin{equation}
\label{troisdemi}
T^{(v_{2}(\mu_{j,w} + 1) )} (\mu_{j,w})=
-1+(\mu_{j,w} + 1)
\left(\frac{3}{2}\right)^{v_{2}(\mu_{j,w} +1)}
,
\end{equation}
and 
the number of waves of $\mathcal{T}$
is infinite,
\item[(iii)] if there exists
$M \geq 2$ such that
any wave
$(\mu_{j,w}=T^{(i_j)}(n),
T^{(i_j + 1)}(n),\ldots,
T^{(i_j +k)}(n))$, $j \geq 1$
and $k \geq 2$,
has a height $k$
which satisfies: $k \leq M$,
the set of prime numbers generated 
by the
wave-minima is infinite, i.e.
$$\bigcup_{j=1}^{+\infty} \mathcal{P}rm
(\mu_{j,w}) 
\hspace{0.5cm}
\mbox{is infinite}.$$

\end{enumerate}
\end{theorem}

\

\section{Proof of the Theorem}
\label{S2}

(i) If
$T^{(j)}(n)$, for some 
$j \geq 0$, is a local minimum of
$\mathcal{T}$
such that
$(T^{(j)}(n), T^{(j+1)}(n))$ is a ripple, we have: 
$T^{(j+1)}(n) > T^{(j)}(n)$ 
and
$T^{(j+2)}(n)=T^{(j+1)}(n)/2$ 
less than $T^{(j)}(n)$.
Each time a ripple-minimum occurs, 
a descent of the trajectory
$\mathcal{T}$ ``occurs" since
$T^{(j+2)}(n) < T^{(j)}(n)$.

(ii) If $\mathcal{T}$
only contains 
ripples and blocks of sequences of even 
integers
successively obtained by division by 2, it 
could not be divergent. Then, we have 
at least one
wave, of a certain height $k$
which is necessarily $\geq 2$, 
in the trajectory $\mathcal{T}$.
The height $k$ is necessarily finite
by the following Proposition; it cannot be infinite.

\begin{proposition}
\label{eqDIOPHprop}
Let $m \geq 3$ be an odd integer.
Let $s \geq 1$ be an integer.
Assume that $T^{(r)}(m)$
is odd for all $r=0, 1,\ldots, s$.
Then
\begin{equation}
\label{eqDIOPH}
3^r (m+1) - 2^r ~=~ T^{(r)}(m) \times 2^r 
\qquad \mbox{for}~~ r=0,1,2, \ldots, s+1.
\end{equation}
\end{proposition}
\begin{proof}
If $r=0$, we have: $(m+1)-1=
T^{(0)}(m)=
m$. Let us prove
\eqref{eqDIOPH} by induction. 
Assume that \eqref{eqDIOPH} is true 
for some $r$ such that $0 \leq r \leq s$.
%with \eqref{eqDIOPH} true, and $j < s$.
Then $T^{(r+1)}(m)$ is given by
$$T^{(r+1)}(m) =
T(T^{(r)}(m))= \frac{3 T^{(r)}(m)+1}{2}
=
\frac{1}{2}\left(
3 \Bigl(
\frac{3^r (m+1)-2^r}{2^r} - 
\Bigr)+1
\right)
.$$
We deduce
$$2^{r+1} T^{(r+1)}(m) =
3^{r+1} (m+1) -3 \times 2^{r} + 2^{r} 
=
3^{r+1} (m+1) - 2^{r+1} .
$$
and the claim.
\end{proof}
Indeed, let us consider a wave-minimum 
$\mu_{j,w}$ of
the trajectory $\mathcal{T}$.
Since $\mu_{j,w}$ is odd, as local minimum,
$\mu_{j,w}+1$ is even. Let
$$ \mu_{j,w}+1 = 2^u p_{1}^{\alpha_1}
p_{2}^{\alpha_2} \ldots p_{q}^{\alpha_q}, \qquad
\mbox{with}~ 
u = v_{2}(\mu_{j,w}+1) \geq 1,$$
denote its prime decomposition
($p_i \neq p_j$ for $i \neq j$ and
$p_i \neq 2$ for $1 \leq i,j \leq q$,
$\alpha_i$s integers).
Since $\mathcal{T}$ is assumed divergent,
all the wave-minima
$\mu_{j,w}$ satisfy:
$\mu_{j,w} \geq 3$.

Let 
$(\mu_{j,w} = T^{(i_j)}(n), T^{(i_j+1)}(n), 
\ldots, T^{(i_j+k)}(n)
= T^{(k)}(\mu_{j,w})
)$
be the wave of height $k$
whose first element is the wave-minimum
$\mu_{j,w}$.
The integers
$T^{(i_j+r)}(n)$, for $r=0$ to $k-1$, are all
odd and $T^{(i_j+k)}(n)$ is even.
Let us show that $k=u$.

Indeed, Proposition \ref{eqDIOPHprop} implies
$$3^r (\mu_{j,w}+1) = (T^{(r)}(\mu_{j,w}) +1)
\times 2^r , \qquad \mbox{for} ~~
0 \leq r \leq k.$$
Hence,
\begin{equation}
\label{wave}
3^r 2^{u-r} p_{1}^{\alpha_1}
p_{2}^{\alpha_2} \ldots p_{q}^{\alpha_q}
= T^{(r)}(\mu_{j,w})+1.
\end{equation}
The prime numbers occurring in the 
prime decomposition of any element
(shifted by $+1$)
of the wave $T^{(r)}(\mu_{j,w})+1$,
$0 < r < k$,
are exactly those occurring in the 
prime decomposition of 
the integer
$\mu_{j,w} +1$, except possibly $3$.
Since $T^{(i_j+k)}(n)$ is even, 
$T^{(i_j+k)}(n)+1$ is odd. We deduce
$2^{u-k}=1$ on the left-hand side
of \eqref{wave}, i.e.
$k = v_{2}(\mu_{1,w}+1)$.

Therefore, from \eqref{wave},
$$T^{(v_{2}(\mu_{j,w} + 1) )} (\mu_{j,w})=
\left(\frac{3}{2}\right)^{v_{2}(\mu_{j,w}+1)}
(\mu_{j,w} + 1)
-1.$$

Now, since the trajectory $\mathcal{T}$ is assumed
divergent, the number of waves of
$\mathcal{T}$, 
each of them of finite height, is necessarily
infinite to have divergence.
The sequence $(\mu_{j,w})_{j \geq 1}$ is infinite.

\

(iii) For proving (iii)
we need to recall a recent theorem of
Bennett and Billerey
(Theorem 3.3 in \cite{bennettbillerey}),
with some notations, 
and adapt it
to the $3x+1$ problem. 
The statement is the following.
  
If $S= \{p_1, p_2 , \ldots, p_q\}$ is a finite 
set of primes, we define the set of
$S$-units to be those rational integers of the 
type  $\pm p_{1}^{\alpha_1} p_{2}^{\alpha_2} 
\ldots p_{q}^{\alpha_q}$, where the $\alpha_i$s
are integers, and $p_i \neq p_j$ for $i \neq j$. 
Let $a$ be a positive integer.

\begin{theorem}[Bennett Billerey 
\cite{bennettbillerey}]
\label{thmBennettBillerey}
There are only finitely
many integer
$S$-units $x$, $y$, $w$
with
$\gcd(x,y) \leq a$ for which there
exists integers $k \geq 2$
and $z \neq 0$ such that
\begin{equation}
\label{twounits}
x+y = w z^k.
\end{equation}
\end{theorem}

In the divergent trajectory
$\mathcal{T}$ all the iterates of
a wave-minima 
$\mu_{j,w}$ in the wave
$(\mu_{j,w} = T^{(i_j)}(n), T^{(i_j+1)}(n), 
\ldots, T^{(i_j+k)}(n)
= T^{(k)}(\mu_{j,w}))$
are linked by
the same $S$-units equation, which
is of the type
\eqref{twounits}, that is
\begin{equation}
\label{eqDIOPHbis}
3^r (\mu_{j,w}+1) - 2^r ~=~ T^{(r)}(\mu_{j,w}) 
\times 2^r 
\qquad \mbox{for}~~ r=0,1,2, \ldots,~ 
k = v_{2}(\mu_{j,w}+1),
\end{equation}
for which
$$x= 3^r (\mu_{j,w}+1),\qquad
y=-2^r,$$
$$w = T^{(r)}(\mu_{j,w}), \qquad z=2, 
\qquad k \geq 2.$$
We have: $\gcd(x,y)= 2^r$. Here, the set $S$
is equal to $\mathcal{P}rm(\mu_{j,w})$ 
and depends upon the
index $j$ of the wave-minimum $\mu_{j,w}$.
The triples
$(3^r (\mu_{j,w}+1), ~-2^r,~
T^{(r)}(\mu_{j,w}))$, 
$0 \leq r \leq v_{2}(\mu_{j,w}+1)$,
are solutions of the finite set
of the $S$-units $x,y,w$ with
$\gcd(x,y) \leq 2^{v_{2}(\mu_{j,w}+1)}$ 
for which there exists
$k \geq 2$ such that
$$x+y=w ~2^k.$$

Let us assume that any wave has a height bounded
from above by an integer $M$.
If $X$ is a finite set of prime numbers,
let us denote by
$U_X$ the set of the $X$-units $x,y,w$ with
$\gcd(x,y) \leq 2^M$ for which there exists
$k \geq 2$ such that
$$x+y=w ~2^k.$$
By Theorem \ref{thmBennettBillerey} the set
$U_X$ is finite.

To prove (iii)  we assume that
$\bigcup_{j=1}^{+\infty} \mathcal{P}rm
(\mu_{j,w})$ is finite. Then we show that this assumption leads to a contradiction.

Let us take:
$$X:= \bigcup_{j=1}^{+\infty} \mathcal{P}rm
(\mu_{j,w}).$$
The set $X$ is finite, and for any $j \geq 1$,
$S:=\mathcal{P}rm(\mu_{j,w})$  is such that
$$S \subset X\quad \mbox{implies}\quad
U_S \subset U_X.$$

Therefore, all the triples 
$(3^r (\mu_{j,w}+1), ~-2^r,~
T^{(r)}(\mu_{j,w}))$, $j \geq 1$,
$0 \leq r \leq v_{2}(\mu_{j,w}+1)$,
belong to $U_X$.
By (ii), it is an infinite collection of triples 
of integers in the finite set $U_X$.
We deduce that there exists
integers $j_1$ and $j_2$,
$j_1 \neq j_2$, and
$r_1 , r_2 \leq M$, such that
$$(3^{r_1} (\mu_{j_1 ,w}+1), ~-2^{r_1},~
T^{(r_1)}(\mu_{j_1 ,w})) =
(3^{r_2} (\mu_{j_2 ,w}+1), ~-2^{r_2},~
T^{(r_2)}(\mu_{j_2 ,w})).$$
We deduce: $r_1 = r_2$.
Considering the first coordinates,
we would have
$\mu_{j_1 ,w}=\mu_{j_2 ,w}$.
This equality means that the trajectory
$\mathcal{T}$ would possess a cycle at infinity,
which would be different from the trivial cycle
$\{1,2\}$. The trajectory 
$\mathcal{T}$ would be bounded.
But the trajectory 
$\mathcal{T}$
is assumed divergent.
Contradiction. We deduce (iii).

\begin{remark}
The form of the equation 
\eqref{eqDIOPHbis} is common 
to all the waves. The only difference from one wave
to another one comes from 
the set $S$ of prime numbers, on which the 
integers $x,y, w$ are decomposed.
\end{remark}

\section{Density and Open Problems}
\label{S3}

With the notations introduced above, 
let us show that the density of the divergent
trajectory 
$\mathcal{T}
=
(T^{(j)}(n))_{j \geq 0}$, for $n \geq 3$ odd, 
(if any) can be defined
using the sets of the prime numbers
generated by the wave-minima
of $\mathcal{T}$.

It can be done as follows: denote
$$I_j:= ~\max {\mathcal Prm}
(\mu_{j,w}),
\qquad j\geq 1.$$
Then we define:

\begin{equation}
\label{density}
\delta_{\mathcal{T}}
:= \limsup_{X \to \infty}~~
\frac{1}{\pi(X)}~ 
\# \Bigl(\bigcup_{1 \leq j, ~I_j \leq X} {\mathcal Prm}
(\mu_{j,w})\Bigr)
\qquad \in[0,1].
\end{equation}

This definition is extended to
any finite trajectory, or to any trajectory ending 
in nontrivial cycles (cf \cite{eliahou}
\cite{eliahou2}): 
in both cases,
$\delta_{\mathcal{T}}=0$ since the 
sets of prime numbers generated by
the wave-minima are finite.
We call $\delta_{\mathcal{T}}$ the density of
$\mathcal{T}$.

As a first step, in the context of
rarefied prime numbers generated by the wave-minima,
a reasonable conjecture 
towards a partial proof of 
the Conjecture of Collatz would be 
the following.

\begin{Conjecture}
Let $\mathcal{T}
=(T^{(j)}(n))_{j \geq 0}$
be the trajectory of the 
odd integer $n \geq 3$
under the Collatz map $T$. Then
$$\delta_{\mathcal{T}} = 0\quad
\Longrightarrow \quad 
\mathcal{T} ~~\mbox{is finite or ends in a nontrivial cycle}.$$
\end{Conjecture}

The second step, towards a proof of the 
Conjecture of Collatz, would be to prove that
the set of trajectories 
$\mathcal{T}$
having a density $\delta_{
\mathcal{T}}>0$ is empty.

\

The combinatorics of
trajectories under $T$ in rational base $3/2$  
has been shown to be simple in
\cite{eliahouvergergaugry}.
It is particularly the case
for strings of successive  odd iterates in a 
given trajectory.
Therefore,
Theorem \ref{mainthm} may justify the need to study
the prime numbers in rational base $3/2$, 
generated by
the local minima along a trajectory.

\vspace{0.5cm}

As a consequence of Theorem  \ref{mainthm},
considering trajectories which possess
waves of arbitrarily large height, we 
formulate the following conjecture.

\begin{Conjecture}
The set of divergent trajectories having
the set of prime numbers generated 
by the wave-minima
$$\bigcup_{j=1}^{+\infty} \mathcal{P}rm
(\mu_{j,w}) $$
finite is empty.
\end{Conjecture}

Acknowledgements.-- The author thanks 
G\"okhan Soydan for useful discussions 
on Diophantine exponential equations and
Shalom Eliahou for his
comments.

\end{document}